\documentclass[noelectric]{jflhack}

\title{A 5-quantifier $(\in,=)$-expression\\ ZF-equivalent to\\ the Axiom of Choice}
\author{Kurt Maes}
\address{Kurt Maes\\ Ertbrug 2\\ 9260 Serskamp\\ Belgium}
\email{maes.kurt@pandora.be}
\primarydata{03E25}
\secondarydata{03E30}
\copyrightinfo{2003} \issue{?}{?}{????}
\received{} \accepted{} \printed{}

\begin{document}
  \begin{abstract}
    In this paper I present an $(\in, =)$-sentence, $AC^{**}$, with only 5 quantifiers, that logically implies the axiom of choice, $AC$.
    Furthermore, using a weak fragment of ZF set theory, I prove that $AC$ implies $AC^{**}$.

    Up to now 6 quantifiers were the minimum and 3 quantifiers don't suffice since all 3-quantifier $(\in, =)$-sentences are decided in a weak fragment of ZF set theory.
    Thus the gap is reduced to the undecided case of a 4 quantifier sentence ZF-equivalent to $AC$.
  \end{abstract}
  \maketitle

  \section{Acknowledgment and Introduction}
    First and foremost I would like to thank Norman Megill for reviewing this paper and its drafts and for checking the proofs. 
    Secondly I would like to thank Harvey Friedman for the challenges he poses to people like myself by posting bold enough conjectures on the
    FOM\endnote{More information can be found at\\ \url{http://www.cs.nyu.edu/mailman/listinfo/fom}.}, an E-mail listing for people interested in the foundations of mathematics.

    In one of his FOM-postings\endnote{
      The posting of Harvey Friedman I'm referring to is entitled ``196:Quantifier complexity in set theory'' and can be found at
      \\\url{http://www.cs.nyu.edu/pipermail/fom/2003-November/007680.html}}
    Harvey Friedman sums up about what was known, up to then, about quantifier complexity in set theory.
    He also makes a number of conjectures.
    I quote the part of interest to this paper below.
    \begin{quote}
      As I said many times on the FOM, all 3 quantifier sentences are decided in a weak fragment of ZF.
      There is a 5 quantifier sentence that is not decided in ZFC, and is provably equivalent to the existence of a subtle cardinal over ZFC.
      \\\ldots\\
      I conjecture that the axiom of choice cannot be stated with 5 quantifiers, but this isn't known even for 4 quantifiers.
      We know that the AxC can be stated with 6 quantifiers (posting \#195).
    \end{quote}
    In this quote a sentence is understood to be part of the primitive language of set theory, which is standard first order predicate calculus with equality and one binary relation, $\in$.
    Furthermore, in counting the number of quantifiers, one counts every individual quantifier, not just alternations of quantifiers.

    This inspired me to investigate the 4- and 5-quantifier sentences, and in particular I set out to disprove the conjecture that there isn't a 5-quantifier $(\in, =)$-sentence equivalent to the axiom of choice.
    I succeeded and this paper presents the proof.
    A draft version of this result was checked by Norman Megill with Metamath, a system for computer-aided formal proof checking developed by Norman Megill.

    All this poses the next question: ``Can $AC$ be stated by an $(\in, =)$-sentence with only 4 quantifiers?''.
    I personally believe this isn't possible.
    This believe comes from the fact that all known (to me) sentences ZF-equivalent to $AC$ have a form like $\forall x \ldots \exists y \phi(x,y)$.
    The dots allow for some premises which must be satisfied by $x$ and $\phi(x,y)$ must, in some way, express that $y$ represents a maximal decision relative to $x$.
    Usually this $\phi$ can't be stated with 2 or less quantifiers.
    The only exception I am aware of is a reworked version of Zorn's Lemma.
    But in this case the premises on $x$ become quantifier-loaded.

    In section 2, I introduce the notion of choice-sets and prove that ``$y$ is a choice-set for $x$'' can be stated with a well-formed-formula (wff) in the primitive language of set-theory containing only 3 quantifiers.
    In particular this wff belongs to the complexity class $\forall\exists\forall$.

    In section 3, I state the axiom of choice, $AC$, in terms of choice-sets.
    Furthermore, I introduce a new ``stronger'' statement $AC^*$, which is stronger in the sense that only first order predicate calculus is required to prove that $AC^*$ implies $AC$.
    This is possible since $AC^*$ exhibits a weaker hypothesis and a stronger conclusion then $AC$.
    Following this, I show that this new statement is implied by $AC$ in a weak fragment of ZF.
    Namely the following:
    \begin{description}
    \item[Extensionality] $\forall x, y (\forall z (z \in x \leftrightarrow z \in y) \rightarrow x = y)$
    \item[Pairing] $\forall x, y \exists z ~z = \lbrace x, y \rbrace$.\\
      In $(\in, =)$-notation this becomes $\forall x, y \exists z \forall a (a \in z \leftrightarrow a = x \vee a = y)$.
    \item[Bounded Separation] If $\varphi$ is a wff in which all quantifiers are bounded, then we have $\forall x \exists y ~y = \lbrace z \in x \mid \varphi \rbrace$.\\
      In $(\in, =)$-notation this becomes $\forall x \exists y \forall z (z \in y \leftrightarrow z \in x \wedge \varphi)$.
    \item[Bounded Replacement] If $\phi$ is a wff in which all quantifiers are bounded and $y$ is not free, then we have $\forall x (\forall z \exists ! z' \phi \rightarrow \exists y \forall z \in x \exists z' \in y ~\phi)$.\\
      In $(\in, =)$-notation this becomes\\
      $\forall x (\forall z \exists z' (\phi \wedge \forall z'' (\phi(z''/z') \rightarrow z'' = z')) \rightarrow \exists y \forall z \exists z' (z \in x \rightarrow z' \in y \wedge \phi))$.
    \end{description}
    The proof is done by describing a method for constructing a set $x^*$ from any set $x$, such that the following properties are satisfied:
    \begin{enumerate}
    \item If $x$ satisfies the hypotheses of $AC^*$, then $x^*$ satisfies the hypotheses of $AC$.
    \item If $y$ is a choice-set for $x^*$, then we can construct a choice-set from this set for $x$ which is not contained in $x$.
    \end{enumerate}

    In section 4, I rewrite $AC^*$ to obtain $AC^{**}$, which is a 5-quantifier $(\in, =)$-sentence.
    The proof of their equivalence only involves first order predicate calculus without any of the axioms of ZF set theory.

  \section{Choice-sets}

  \begin{definition}
    We call a set $y$ a choice-set for a set $x$, and write $C(y,x)$, if for any non-empty element, $z \in x$, the intersection of $z$ and $y$, $z \cap y$, is a singleton.
    I.e. $y$ chooses exactly one element from every non-empty $z$ in $x$.
    Formally this could be written as
    $$\forall z \in x (z \neq \emptyset \rightarrow \exists ! a \in y ~a \in z) \text{.}$$
    If we want to use $(\in, =)$-notation, this can be translated into
    $$\forall z \in x (\exists b ~b \in z \rightarrow \exists a (a \in z \wedge a \in y \wedge \forall b (b \in z \wedge b \in y \rightarrow b = a))) \text{.}$$
  \end{definition}

  \begin{lemma}
    Suppose we have three wff's $X(t)$, $Y(t)$ and $Z(r,t)$, in which the variables $a$ and $b$ do not occur freely.
    Suppose furthermore, that $\forall t (Y(t) \rightarrow X(t))$ is valid.
    If we now define the wff's $A$ and $B$ as below, then $A$ and $B$ are equivalent.

    \begin{itemize}
    \item $A \equiv \exists b X(b) \rightarrow \exists a (Y(a) \wedge \forall b Z(b,a))$,
    \item $B \equiv \exists a \forall b [(X(b) \rightarrow Y(a))\wedge (Y(a)\rightarrow Z(b,a))]$.
    \end{itemize}
  \end{lemma}
  \begin{proof}~\\
    We need to prove $A\vdash B$ and $B\vdash A$

    \begin{description}
    \item[$A\vdash B$]
      Suppose that we have $A$.
      We have to consider two possibilities.
      \begin{description}
      \item[$\exists b X(b)$]
        In this case, $A$ gives us an $a$ that satisfies both $Y(a)$ and $\forall b Z(b,a)$.
        Now take any $b$.
        Since we have $Y(a)$, we definitely have $X(b) \rightarrow Y(a)$.
        Because of $\forall b Z(b,a)$, we also have $Z(b,a)$, which in turn implies $Y(a) \rightarrow Z(b,a)$.
        Hence for our chosen $a$ and any $b$, we find that we have $(X(b) \rightarrow Y(a)) \wedge (Y(a) \rightarrow Z(b,a))$.
      \item[$\neg \exists b X(b)$]
        In this case we have $\forall b \neg X(b)$.
        Now take any $a$ and $b$.
        By assumption, we have $\neg X(b)$, which gives us $X(b) \rightarrow Y(a)$.
        Also by assumption, we have $\neg X(a)$.
        Since we required $\forall t (Y(t) \rightarrow X(t))$ to be valid, we have $Y(a) \rightarrow X(a)$.
        These two together give us $\neg Y(a)$, which result in the validity of $Y(a) \rightarrow Z(b,a)$.
        Again $B$ follows.
        We even got $B$ with a universal quantifier instead of an existential quantifier.
      \end{description}
    \item[$B\vdash A$]
      Suppose we now have $B$.
      To prove $A$, we assume $\exists b X(b)$.
      Hence we may assume a $b$ that satisfies $X(b)$.
      $B$ gives us an $a$, that satisfies $\forall b [(X(b) \rightarrow Y(a))\wedge (Y(a)\rightarrow Z(b,a))]$.
      So in particular our chosen $a$ and $b$ from above satisfy both $X(b)$ and $X(b) \rightarrow Y(a)$, from which we can conclude $Y(a)$.
      Now take any $b$.
      Again because of our choice for $a$, we get that $a$ and $b$ satisfy $Y(a)\rightarrow Z(b,a)$.
      But since we have $Y(a)$, this implies $Z(b,a)$.
      Since $b$ was chosen arbitrarily, our chosen $a$ also satisfies $\forall b Z(b,a)$.
      This concludes the proof.
    \end{description}
  \end{proof}

  \begin{corollary}
    Notice that that the premises of the above lemma are fulfilled when we choose $X(t)$, $Y(t)$ and $Z(r,t)$ as below.
    \begin{itemize}
    \item $X(t) \equiv t \in z$,
    \item $Y(t) \equiv X(t) \wedge t \in y$,
    \item $Z(r,t) \equiv Y(r) \rightarrow r = t$.
    \end{itemize}
    This gives the following $A$ and $B$:
    \begin{itemize}
    \item $A \equiv \exists b ~b \in z \rightarrow \exists a (a \in z \wedge a \in y \wedge \forall b (b \in z \wedge b \in y \rightarrow b = a))$,
    \item
      \begin{tabular}{rcl}
        $B \equiv \exists a \forall b [$ &          & $(b \in z \rightarrow a \in z \wedge a \in y)$                                     \\
                                         & $\wedge$ & $(a \in z \wedge a \in y \rightarrow (b \in z \wedge b \in y \rightarrow b = a))]$ \\
      \end{tabular}.\\
    \end{itemize}
    $A$ clearly states: If $z$ is non-empty, then the intersection between $z$ and $y$ is a singleton.
    While $B$ shows that this can be said with only two quantifiers.
    More in particular we may conclude that $C(y,x)$ is equivalent to
    $$\forall z (z \in x \rightarrow B(x,y,z)) \text{.}$$
    Hence ``$y$ is a choice-set for $x$'' can be stated in $(\in, =)$-notation with as little as 3 quantifiers.
  \end{corollary}

  \section{Axiom of Choice}

  \begin{definition}
    One of the formulations of the axiom of choice, $AC$, states:\\
    \emph{For any set $x$ consisting of non-empty pairwise disjoint elements, there exists a choice-set $y$ for $x$.}\\
    Formally this can be written as:
    $$AC \equiv \forall x [AC_{h,1}(x) \wedge AC_{h,2}(x) \rightarrow \exists y C(y,x)] \text{,}$$
    where we have
    $$AC_{h,1}(x) \equiv \forall z \in x ~z \neq \emptyset \text{ and}$$
    $$AC_{h,2}(x) \equiv \forall z, z' \in x (z \neq z' \rightarrow z \cap z' = \emptyset) \text{.}$$

    In the literature, one might find an alternative formulation for the axiom of choice that does not require the elements of $x$ to be non-empty to guarantee the existence of a choice-set.
    In the presence of the Axiom of Bounded Separation, these two formulations are equivalent, since $y$ is a choice-set for $x$ iff $y$ is a choice-set for $\lbrace z \in x \mid z \neq \emptyset \rbrace = \lbrace z \in x \mid \exists a \in z ~a = a \rbrace$.
    This follows easily from the fact that $C(y,x)$ only gives information about the non-empty elements of $x$.
  \end{definition}

  \begin{definition}
    In what follows, we will be interested in the following statement:
    $$AC^* \equiv \forall x (AC_h^*(x) \rightarrow \exists y (y \notin x \wedge C(y,x))) \text{,}$$
    where we have
    $$AC_h^*(x) \equiv \forall z \in x \exists a \in z \forall z' \in x (z \neq z' \rightarrow a \notin z') \text{.}$$
    I.e. $AC_h^*(x)$ states that all elements of $x$ contain an element not contained in any other element of $x$.
    The  main purpose of this section is to prove the equivalence of this statement with $AC$.
  \end{definition}

  \begin{lemma}
    $AC_{h,1}(x) \wedge AC_{h,2}(x) \rightarrow AC_h^*(x)$.
  \end{lemma}
  \begin{proof}
    Take any $z \in x$.
    Since we have $AC_{h,1}(x)$, we have an $a$ in $z$.
    Now, $AC_{h,2}(x)$ gives us that this $a$ is not an element of
    any other element of $x$.
    Hence $AC_h^*(x)$ follows.
  \end{proof}

  \begin{theorem}
    $AC^* \rightarrow AC$.
  \end{theorem}
  \begin{proof}
    Suppose $AC^*$ is valid.
    Now suppose $x$ satisfies the hypotheses of $AC$.
    The previous lemma states that $x$ then also satisfies the hypotheses of $AC^*$.
    Hence by $AC^*$, we have a choice-set $y$ for $x$, which even is guaranteed not to be an element of $x$.
  \end{proof}

  Notice that the proof of $AC^* \rightarrow AC$ did not use any of the axioms of set-theory.
  The reverse implication, on the other hand, requires some of the other axioms of set-theory.

  In what follows, I will indicate in parentheses which axioms and which previous results are required to prove the stated result.

  \begin{lemma}{(Extensionality, Bounded Separation)}
    Given any set $x$, then $\phi(z,z_x)$ defined below is a function in the sense of the Axiom of Bounded Replacement.
    $$\phi(z,z_x) \equiv z_x = \lbrace a \in z \mid \forall z^* \in x (z \neq z^* \rightarrow a \notin z^*) \rbrace$$
    Furthermore, for this function, we have the following properties:
    \begin{enumerate}
    \item $\forall z \in x ~z_x \subseteq z$,
    \item $\forall z, z' \in x (z \neq z' \rightarrow z_x \cap z' = \emptyset)$,
    \item $\forall z, z' \in x (z \neq z' \rightarrow z_x \cap z'_x = \emptyset)$,
    \item $AC_h^*(x) \rightarrow \forall z \in x ~z_x \neq \emptyset$.
    \end{enumerate}
    Here $z_x$ (resp.~$z'_x$) denotes the image of $z$ (resp.~$z'$) for the function $\phi$ defined above.
    Formally this means that property 3 actually stands for
    $\forall z, z' \in x \forall z_x, z'_x [\phi(z,z_x) \wedge \phi(z',z'_x) \rightarrow (z \neq z' \rightarrow z_x \cap z'_x = \emptyset)]$
  \end{lemma}
  \begin{proof}
    Because of the Axiom of Bounded Separation we have for all $z$, a set $z_x$ such that $\phi(z,z_x)$ is satisfied.
    The Axiom of Extensionality guarantees the uniqueness of this $z_x$.
    Hence $\phi$ does indeed define a function.
    Since $\phi(z,z_x)$ can be stated in $(\in, =)$-notation with only bounded quantifiers as follows:
    \begin{eqnarray}
      \nonumber \phi(z,z_x) \equiv &      & \forall a \in z (\forall z^* \in x (z \neq z^* \rightarrow a \notin z^*) \rightarrow a \in z_x) \\
      \nonumber                    &\wedge& \forall a \in z_x (a \in z \wedge \forall z^* \in x (z \neq z^* \rightarrow a \notin z^*)) \text{,}
    \end{eqnarray}
    $\phi$ does indeed satisfy all the conditions of the Axiom of Bounded Replacement.\\
    We now prove the properties of this function.
    \begin{enumerate}
    \item This follows immediately from the definition of $z_x$.
    \item Suppose this property didn't hold. Hence we have two distinct elements, $z$ and $z'$, in $x$ and an $a$, such that $a \in z_x \cap z$ holds.
      By the previous property, we would have $a \in z \cap z'$.
      But by definition of $z_x$, $a$ would then not be an elements of $z_x$.
      Which is in contradiction with our chosen $a$.
    \item Given any two distinct elements, $z$ and $z'$, in $x$, the two previous properties allow us to derive $z_x \cap z'_x \subseteq z_x \cap z' = \emptyset$.
    \item Suppose $x$ satisfies $AC_h^*(x)$.
      Now take any $z \in x$.
      $AC_h^*(x)$ gives us an $a \in z$, not contained in any $z'$ in $x$ different from $z$.
      Hence by definition we have $a \in z_x$.
    \end{enumerate}
  \end{proof}

  \begin{corollary}(Bounded Replacement, Previous Lemma)
    The Axiom of Bounded Replacement and the previous lemma guarantee the existence of a set $x^*$ as the image of the function $\phi$ restricted to $x$.
    $$x^* = \lbrace z_x \mid z \in x \rbrace$$
    One easily verifies that such an $x^*$ satisfies the following properties:
    \begin{itemize}
    \item $AC_{h,2}(x^*)$ (corresponds to property 3 in the previous lemma),
    \item $AC_h^*(x) \leftrightarrow AC_{h,1}(x^*)$ (``$\rightarrow$'' corresponds to property 4 in the previous lemma, and the reverse follows from properties 1 and 2 in the previous lemma).
    \end{itemize}
  \end{corollary}

  \begin{theorem}(Pairing, Bounded Separation, Previous corollary, Previous lemma)
    $AC \rightarrow AC^*$.
  \end{theorem}
  \begin{proof}
    Suppose $AC$ is valid and $x$ is a set that satisfies $AC_h^*(x)$.
    Now consider a set $x^*$ given by the above corollary.
    This same corollary and $AC_h^*(x)$ guarantees that $x^*$ satisfies the hypotheses of $AC$.
    Hence $AC$ gives us a choice-set $y$ for $x^*$.
    The Axiom of Bounded Separation guarantees us the existence of a set
    $$y' = \lbrace a \in y \mid \exists z \in x ~a \in z_x \rbrace \text{.}$$
    One easily verifies that we have
    $$y' = y \cap (\bigcup x^*) \text{.}$$
    On the other hand, the previous lemma allows one to verify that
    $$z_x = z \cap (\bigcup x^*)$$
    is valid for any $z \in x$.
    This allows us to derive
    $$z \cap y' = z \cap (y \cap (\bigcup x^*)) = (z \cap (\bigcup x^*)) \cap y = z_x \cap y \text{.}$$
    Since $y$ is a choice-set for $x^*$, we find that $z \cap y'$ is a singleton for all $z \in x$.
    Hence $y'$ is a choice-set for $x$.

    If $y'$ is not an element of $x$, then nothing remains to be proven.
    So suppose $y'$ is an element of $x$.
    Since $y'$ is a choice-set for $x$, we find that $y' = y' \cap y' = \lbrace a \rbrace$ for some $a$.
    However, since we have $AC_h^*(x)$, $a$ cannot be contained in any other element of $x$ and no element of $x$ is empty.
    On the other hand, since $y'$ is a choice-set for $x$, $a$ would have to be contained in all non-empty elements of $x$.
    Hence we get $x = \lbrace y' \rbrace$.
    Now if we where to find a set $b$ different from $a$, then the Axiom of Pairing gives a set $y'' = \lbrace a, b \rbrace$.
    Such a set would then still be a choice-set for $x$ and it would not be an element of $x$.

    The search for this set $b$ can be done by different means (read: using different axioms of set-theory).
    Since we have already been using the Axiom of Bounded Separation, we will use this route.
    If $a$ is empty, then $y'$ which is non-empty will do.
    If on the other hand $a$ is non-empty, then the Axiom of Bounded Subsets ($\lbrace u \in a \mid u \neq u \rbrace$) guarantees that there exists an empty set $b$ as a subset of $a$.
    This $b$ will do in this case.
  \end{proof}

  \begin{corollary}
    If we summarize, then we find that to prove the equivalence of $AC$ and $AC^*$ we used the following axioms of set-theory.
    \begin{itemize}
    \item Axiom of Extensionality,
    \item Axiom of Pairing,
    \item Axiom of Bounded Separation,
    \item Axiom of Bounded Replacement.
    \end{itemize}

    Since $x^*$ is a subset of the powerset of $x$, one could also replace the Axiom of Bounded Replacement by the Axiom of Powersets.

    Note that while we did use $\cap$, $\cup$ and $\subseteq$ quite a bit, these were always used to construct classes and their relationships as classes.
    Hence their use can be eliminated from the above, be it at the cost of readability.
  \end{corollary}

  \section{5-quantifier Axiom of Choice}

  The following result requires only first-order predicate calculus, i.e. none of the axioms of ZF were used.

  \begin{theorem}
    The sentence
    $$AC^{**} \equiv \forall x \exists y \forall z \exists a \forall b [(y \in x \wedge A(x,y,z,a)) \vee (y \notin x \wedge B(x,y,z,a,b))]$$
    where $A(x,y,z,a)$ and $B(x,y,z,a,b)$ are given by
    $$A(x,y,z,a) \equiv z \in y \rightarrow a \in x \wedge a \neq y \wedge z \in a$$
    and
    \begin{eqnarray}
      \nonumber B(x,y,z,a,b) \equiv z \in x \rightarrow &      & (b \in z \rightarrow a \in z \wedge a \in y)\\
      \nonumber                                         &\wedge& (a \in z \wedge a \in y \rightarrow (b \in z \wedge b \in y \rightarrow b=a))
    \end{eqnarray}
    is equivalent to $AC^*$.\\
  \end{theorem}
  \begin{proof}~\\
    \begin{enumerate}
    \item
      By corollary 2.3, we find that
      $$\exists y (y \notin x \wedge C(y,x))$$
      is equivalent to
      \begin{eqnarray}
        \nonumber \exists y (y \notin x \wedge \forall z (z \in x \rightarrow \exists a \forall b [ &      & (b \in z \rightarrow a \in z \wedge a \in y) \\
        \nonumber                                                                                   &\wedge& (a \in z \wedge a \in y \rightarrow (b \in z \wedge b \in y \rightarrow b = a))]))
      \end{eqnarray}
      Since $z \in x$ doesn't mention $a$ nor $b$, this is equivalent to
      $$\exists y (y \notin x \wedge \forall z \exists a \forall b B(x,y,z,a,b)) \text{.}$$
      Since $y \notin x$ doesn't mention $z$, $a$ nor $b$, this is equivalent to
      $$\exists y \forall z \exists a \forall b (y \notin x \wedge B(x,y,z,a,b)) \text{.}$$
    \item
      $\neg AC_h^*(x)$ can be rewritten as
      $$\neg \forall z [z \in x \rightarrow \exists a (a \in z \wedge \forall z' [z' \in x \rightarrow (z \neq z' \rightarrow a \notin z')])] \text{.}$$
      Carrying the negation through the quantifiers and the logical connections on its way, we find this to be equivalent to
      $$\exists z [z \in x \wedge \forall a (a \in z \rightarrow \exists z' [z' \in x \wedge z \neq z' \wedge a \in z'])] \text{.}$$
      Since $a \in z$ does not mention $z'$, this is equivalent to
      $$\exists z [z \in x \wedge \forall a \exists z' (a \in z \rightarrow [z' \in x \wedge z \neq z' \wedge a \in z']])] \text{.}$$
      Since $z \in x$ does not mention $a$ nor $z'$, this is equivalent to
      $$\exists z \forall a \exists z' [z \in x \wedge (a \in z \rightarrow [z' \in x \wedge z \neq z' \wedge a \in z']])] \text{.}$$
      Replacing the variables $z$, $a$ and $z'$ simultaneously with the variables $y$, $z$ and $a$, this becomes
      $$\exists y \forall z \exists a (y \in x \wedge A(x,y,z,a)) \text{.}$$
      Since this wff does not mention $b$, this is equivalent to
      $$\exists y \forall z \exists a \forall b (y \in x \wedge A(x,y,z,a)) \text{.}$$
    \item
      One easily verifies that $\neg [(y \in x \wedge A(x,y,z,a)) \wedge (y \notin x \wedge B(x,y,z,a,b))]$ is always valid.
    \end{enumerate}
    The above points allow us to verify our theorem.
    This goes like this:
    $$AC^*$$
    $$\updownarrow$$
    $$\forall x (\neg AC_h^*(x) \vee \exists y (y \notin x \wedge C(y,x)))$$
    $$\updownarrow (1,2)$$
    $$\forall x [\exists y \forall z \exists a \forall b(y \in x \wedge A(x,y,z,a)) \vee \exists y \forall z \exists a \forall b (y \notin x \wedge B(x,y,z,a,b))]$$
    $$\updownarrow (3)$$
    $$\forall x \exists y \forall z \exists a \forall b [(y \in x \wedge A(x,y,z,a)) \vee (y \notin x \wedge B(x,y,z,a,b))]$$
  \end{proof}

  \begin{corollary}
    Hence we have an $(\in, =)$-sentence with only 5 quantifiers, $AC^{**}$, equivalent to the axiom of choice.
  \end{corollary}

  \begin{remark}
    There is a way to obtain a shorter 5-quantifier $(\in, =)$-sentence equivalent to the axiom of choice.
    This follows from the following observations.
    \begin{itemize}
    \item If $x$ is a set without empty elements, then $y$ is a choice-set for $x$ iff the intersection of $y$ and any element of $x$ (not restricted to non-empty elements of $x$) is a singleton.
    \item If $x$ satisfies $AC_h^*(x)$, then it does not contain any empty elements.
    \end{itemize}
    Hence we find that $AC^*$ is equivalent to $\overline{AC}^*$ obtained from $AC^*$ by replacing $C(y,x)$ with $\overline{C}(y,x)$ given by:
    $$\overline{C}(y,x) \equiv \forall z \in x \exists a (a \in z \wedge a \in y \wedge \forall b (b \in z \wedge b \in y \rightarrow b = a))) \text{.}$$
    We then could translate $\overline{AC}^*$ into $\overline{AC}^{**}$, as we did translate $AC^*$ into $AC^{**}$, not needing the step using corollary 2.3.
    This would result in replacing $B(x,y,z,a,b)$ with
    $$\overline{B}(x,y,z,a,b) \equiv z \in x \rightarrow a \in z \wedge a \in y \wedge (b \in z \wedge b \in y \rightarrow b=a) \text{.}$$
    This shortens our sentence by 16 symbols (4 parentheses, 3 logical connectors and 3 atomic formulas).
  \end{remark}
\end{document}